\documentclass{article}
\usepackage{sonar-radon,graphicx}
\title{Operator identities relating sonar and Radon transforms in Euclidean space}
\author{Aleksei Beltukov\\
University of the Pacific\\
\texttt{abeltuko@pacific.edu}\\
and\\
David Feldman\\
University of New Hampshire\\
\texttt{David.Feldman@unh.edu}
}
\date{\today}
\begin{document}
\maketitle
\begin{abstract}
  We establish new relations which connect Euclidean sonar
transforms (integrals taken over spheres with centers in a
hyperplane) with classical Radon transforms. The relations, stated
as operator identities, allow us to reduce the inversion of sonar
transforms to classical Radon inversion.

\end{abstract}
\section{Introduction}
\label{s:Introduction}

As we aim to relate sonar transforms with Radon transforms, we
must begin by recalling key definitions which in turn requires us
first to fix some notation.
\begin{math}
  \h{n} := \r{n-1} \times (0,\infty)
\end{math}
denotes the upper half space of $\r{n}$. Points in $\h{n}$ will be
written as $(x,y)$ with
\begin{math}
  x \in \r{n-1}
\end{math}
and $y>0$. We write $|x|$ for the Euclidean vector norm of
\begin{math}
  x \in \r{n-1},
\end{math}
$dt$ for the Euclidean volume element on $\r{n-1}$.
\begin{math}
  \cc{\h{n}}
\end{math}
(resp.
\begin{math}
  \ccc{\h{n}}
\end{math})
denotes the set of smooth compactly supported functions (resp.
smooth functions) on $\r{n}$ supported in $\h{n}$. $\s[x]{y}{n-1}$
will denote the sphere in $\r{n}$ centered at $x$ with radius $y$
(empty if $y<0$) carrying area measure $\dS{}$.
\par
We now define the sonar transform
\begin{math}
  \sonar: \cc{\h{n}} \mapsto \ccc{\h{n}}
\end{math}
as follows.
\begin{definition} \label{d:sonar}
Given
\begin{math}
  f \in \cc{\h{n}}
\end{math},
\begin{equation*}
  \sonar[f](x,y) := \int_{\s[x]{y}{n-1}} f \, \dS{}.
\end{equation*}
The \textit{centerset variable} $x$ parameterizes the
\textit{centerset}
\begin{math}
  \r{n-1} \times \{0\}.
\end{math}
On occasion we call $y$ the \textit{radial variable}.
\end{definition}
More explicitly,
\begin{equation} \label{e:sonar}
  \sonar[f](x,y) = \int_{|t| < y}
    f \left( x + t, \sqrt{y^{2} - |t|^{2}} \right) \,
    \frac{y \, dt}{\sqrt{y^{2} - |t|^{2}}}.
\end{equation}
The sonar data $\sonar[f]$ generally does not have compact
support.  However the restriction of $\sonar[f]$ to any hyperplane
parallel with the centerset is compactly supported which justifies
various compositions of transforms below.
\begin{remark} \label{r:support}
While we restrict $\sonar$ to $\cc{\h{n}}$ for the sake of our
subsequent derivations, Definition \ref{d:sonar} makes sense for
locally integrable functions $f$.
\end{remark}
Courant and Hilbert initiated the study of  $\sonar$  in
\textit{Methods of Mathematical Physics, Volume II}
\cite{Courant-Hilbert-1989}, where they established its
injectivity on the space of continuous functions and used the
result to investigate hyperbolic partial differential equations.
While \cite{Courant-Hilbert-1989} terms the mapping $\sonar$
``integrals over spheres centered in the plane,'' our more
efficient ``sonar'' terminology follows recent applications of
$\sonar$ to \textit{marine tomography} as in work of Louis and
Quinto \cite{Louis-Quinto-2000}, where the operator $\sonar$ (in
dimension three) models naval sonar data.
\par
Operator $\sonar$  has other practical uses. As Cheney
\cite{Cheney-2001} explains, in dimension two $\sonar$ describes
synthetic aperture radar.  As $\sonar$ and its generalizations
abstract the behavior of reflected waves (echoes) whether
acoustic, electromagnetic, or mechanical, they play a central role
in \textit{reflective tomography}, including marine tomography and
radar theory.
\par
For the sake of  recalling the classical Radon  transform
$\radon$, we will write the set of all hyperplanes in a real
vector space $V$ as ${\cal P}(V)$, which  carries the structure of
a smooth manifold.  We now define
\begin{math}
  \radon_V:\cc{V}\mapsto \cc{{\cal P}(V)}
\end{math}
as follows.
\begin{definition} \label{d:Radon}
\begin{equation*}
  \radon_{V}[f](P) := \int_P f(x) \, \dm{x},
\end{equation*}
where $\dm{x}$ denotes planar surface measure on
\begin{math}
  P \in {\cal P}(V).
\end{math}
\end{definition}
For $V=\r{n}$ we write the Radon transform simply as $\radon$
whereas for $V=\r{n-1} \times \{0\} \subset \r{n}$ we denote the
Radon transform $\radonc$ and call it the \textit{centerset Radon
transform}. Note that both the sonar and Radon transforms reduce
to the identity operator when $n=1$.
\par
The operators $\sonar$ and $\radon$ appear quite different
conceptually. For example, Helgason \cite{Helgason-1980} shows
that the Radon transform  has a very special group-theoretic
structure which leads to the following inversion formula in
$\r{n}$ (taken from \cite{Helgason-1980}, p.15):
\begin{equation} \label{e:iRadon}
  \frac{(4 \, \pi)^{\frac{n-1}{2}} \, \Gamma \left( \frac{n}{2} \right)}
       {\Gamma \left( \frac{1}{2} \right)} \,
  f(x) =
  (-\Delta)^{\frac{n-1}{2}}
  \left( \radon^{*} \circ \radon \right)[f](x)
\end{equation}
Here the superscript $*$ indicates the adjoint operator while
$\Delta$ stands for the Laplacian. In odd dimensions,
\begin{math}
  (-\Delta)^{\frac{n-1}{2}}
\end{math}
is a differential operator; in even dimensions, the fractional
power of the negative Laplacian takes the form of a
pseudodifferential operator that should be interpreted in terms of
Riesz potentials.
\par
In sharp contrast, the composition
\begin{math}
  \sonar^{*} \circ \sonar
\end{math}
does not exist, as $\sonar[f]$ may lack compact support even when
$f$ has one. And, despite the sonar transform's sizable symmetry
group (on $\h{n}$, a semidirect product
\begin{math}
  O(n-1) \propto \r{n-1}
\end{math}
of the orthogonal group with translations), at present we lack a
group-theoretic interpretation of $\sonar$.
\par
In light of these  differences, unexpected close relations between
these two operators carry intrinsic interest. Denisjuk
\cite{Denisjuk-1999}
found the first such relation:
\begin{theorem}[Denisjuk, 1999] \label{t:Denisjuk}
$\b{1}{n}$ denotes the unit ball in $\r{n}$. There exists a
mapping
\begin{math}
  \phi: \b{1}{n} \mapsto \h{n}
\end{math}
(related to  stereographic projection) and a certain non-negative
weight
\begin{math}
  \weight: \b{1}{n} \to \r{}
\end{math}
 such that for any
smooth compactly supported function
\begin{math}
  f \in \cc{\h{n}},
\end{math}
\begin{equation} \label{e:Denisjuk}
  \sonar[f] = \radon[\weight \cdot (f \circ \phi)].
\end{equation}
\end{theorem}
Using Equation \eqref{e:Denisjuk}, Denisjuk expressed the inverse
$\sonar^{-1}$ as a pull-back of the inverse Radon transform
$\radon^{-1}$ and established Plancherel identities for sonar.
Theorem \ref{t:Denisjuk} was later used in \cite{Palamodov-2000}
by Palamodov to pull back various microlocal estimates and perform
$\Lambda$-type reconstruction on integrals over arcs.  As both
\cite{Denisjuk-1999} and \cite{Palamodov-2000} demonstrate,
sonar-Radon relations of type \eqref{e:Denisjuk} can be
effectively used to translate any Radon result into a
corresponding sonar statement.
\par
Denisjuk connects the sonar transform of a given function to the
Radon transform of a \textit{different} function with, \textit{a
priori}, a different support.  Since
\begin{math}
  \h{n}\subset \r{n}
\end{math}
and thus
\begin{math}
  \cc{\h{n}}\subset \cc{\r{n}},
\end{math}
it makes sense to speak of the Radon transform and sonar
transform of one and the same function $f$ (as long as $f$ has
support in $\h{n}$).  This article addresses the following natural
question:
\begin{quote}
  \emph{How can one pass directly from the sonar transform of a function to its Radon transform?}
\end{quote}
\par
A simplified version of our main result says that for almost all
planes $P$, we can compute $\radon[f](P)$ as
\begin{displaymath}
 (\W \circ \radonw[1/y] \circ \D{\frac{n-2}{2}}
  \circ \radonc \circ \sonar)[f](P).
\end{displaymath}
for certain explicitly described operators $\W$, $\radonw[1/y]$,
and $\D{\frac{n-2}{2}}$ each of which has a natural geometric
meaning.

\section{Main Result}
\label{s:MainResult}

We aim to compute the function $\radon[f]$ explicitly from the
function $\sonar[f]$, the \textit{sonar data} associated to $f$.
In order to make our calculations as explicit as possible, we need
a suitable parameterization of ${\cal P}(\r{n})$. Our formula for
calculating $\radon[f](P)$, as it turns out, breaks into various
cases depending on the geometry of the hyperplane $P$ relative to
the centerset of the sonar transform. Accordingly, we write
\begin{displaymath}
  {\cal P}(\r{n})= h \cup v \cup s,
\end{displaymath}
a disjoint union, with
\begin{itemize}
  \item[$h$] the set of planes parallel to $\r{n-1}$;
  \item[$v$] the set of planes perpendicular to $\r{n-1}$; and
  \item[$s$] the set of all other planes (the \textit{slanted} planes).
\end{itemize}
We shall write $\radonh[f]$, $\radonv[f]$, and $\radons[f]$ for
the corresponding restrictions of $\radon[f]$ to $h$, $v$ and $s$.
\par
To make matters more precise, we begin with a general
parameterization of ${\cal P}$.  To a pair $(\omega,p)$,
\begin{math}
  \omega \in \s{1}{n-1},
\end{math}
$p \geq 0$, we may associate the hyperplane
\begin{math}
  P = \{ x \in \r{n} | \omega \cdot x = p \};
\end{math}
here $\cdot$ denotes the standard inner product in $\r{n}$.
$P$~almost determines $(\omega,p)$; only if $p=0$, $\omega$ may
vary by a sign.
\par
Now we adapt this framework to take account of the centerset.  We
write
\begin{math}
  \omega = (\omega',\omega_{n}) \in \r{n-1} \times \r{},
\end{math}
and similarly
\begin{math}
  x = (x',x_{n}).
\end{math}
The equation
\begin{math}
  \omega \cdot x = p
\end{math}
now takes the form
\begin{math}
  \omega' \cdot x' + \omega_{n} \, x_{n} = p.
\end{math}
We can now distinguish three cases, as above:
\begin{itemize}
  \item[$(h)$] $\omega'$ vanishes;
  \item[$(v)$] $\omega_{n}$ vanishes;
  \item[$(s)$] all others.
\end{itemize}
Start with $(h)$.  Since
\begin{math}
  \omega \in \s{1}{n-1},
\end{math}
we must have $\omega_{n}=\pm 1$. Dividing through by $\omega_{n}$,
the defining equation of the plane has the form $x_{n}=y$ (for
some appropriate $y$). The variable $y$ can now parameterize
$(h)$. If
\begin{math}
  f \in \cc{\h{n}}
\end{math}
then
\begin{equation} \label{e:Radonh}
  \radonh[f](y) =
  \int_{x \in \r{n-1}}
    f(x,y) \,
  dx.
\end{equation}
Thus we integrate over the centerset variable.
\par
Now turn to  $(v)$.  A plane in $v$ has defining equation
\begin{math}
  \omega' \cdot x'= p.
\end{math}
The intersection of a vertical plane with the centerset determines
it, so parameterizing $v$ amounts to parameterizing the set of all
hyperplanes in the centerset.  As before, $P$ almost determines
$(\omega',p)$; only if $p=0$, $\omega'$ may vary by a sign.
Accordingly, we shall write
\begin{equation} \label{e:Radonv}
  \radonv[f](\omega,p) =
  \int_{0}^{\infty}
  \left[
    \int_{\omega \cdot x = p}
      f(x,y) \,
    \dm{x}
  \right] \,
  dy
\end{equation}
where $\omega=(\omega',0)$.  Observe that
\begin{equation} \label{e:Radonc}
  \int_{\omega \cdot x = p}
    f(x,y) \,
  \dm{x}=
  \overline{\radon}[f(\cdot,y)](\omega,p)
\end{equation}
(which we view  as an equality between functions of $y$), so we
may also write Equation \eqref{e:Radonv} as
\begin{displaymath}
  \radonv[f](\omega,p) =
  \int_{0}^{\infty}
    \overline{\radon}[f(\cdot,y)](\omega,p) \,
  dy.
\end{displaymath}

The final case $(s)$ comprises a dense open subset of ${\cal
P}(\r{n})$, and thus makes the most substantial contribution to
our union (provided $n>1$).
\par
A plane in $(s)$ has defining equation
\begin{math}
  \omega' \cdot x' + \omega_{n} \, x_{n} = p
\end{math}
with neither $\omega'$ nor $\omega_{n}$ vanishing. We can scale
this equation so as to normalize $\omega'$ and simultaneously
render the coefficient of $x_{n}$ negative. Thus we can
unambiguously choose a defining equation for the same plane with
\begin{displaymath}
  (\omega',p,-\omega_{n}) \in \s{1}{n-2} \times \r{} \times \h{}
\end{displaymath}
Now the  pair $(\omega',p)$ by itself determines the intersection
of $P$ with the centerset, so $\omega_{n}$ controls the angle
between $P$ and the centerset.  More explicitly, write
\begin{math}
  \omega_{n} = -\cot\beta
\end{math} with
\begin{math}
  \beta \in (0,\pi/2).
\end{math}
Then the defining equation of the plane has the form
\begin{displaymath}
  \omega' \cdot x' -\cot\beta \, x_{n} = p.
\end{displaymath}
Intersecting with the parallel translate of the centerset where
\begin{math}
  x_{n} = s \, \sin\beta
\end{math}
gives
\begin{displaymath}
  \omega' \cdot x'= p + s \, \cos\beta.
\end{displaymath}
Henceforth we parameterize $(s)$ by $(\omega',p,\beta)$; $\beta$
now directly represents the angle between $P$ and the centerset.
\par
Integration of
\begin{math}
  f \in \cc{\h{n}}
\end{math}
over $P$ can be split into integration over
\begin{math}
  P \cap \r{n-1}
\end{math}
followed by integration over $\beta$. Explicitly,
\begin{equation} \label{e:Radons}
  \radons[f](\omega',p,\beta) =
  \int_{0}^{\infty}
    \overline\radon[f(\cdot,s\sin \beta)]
    \left( \omega, p + s \, \cos \beta\right)
  ds.
\end{equation}
\begin{remark} \label{r:VnotS}
We avoid including $(v)$ in $(s)$ as a special case $\beta=\pi/2$
both to get a good parameterization of $(s)$ and because the cases
require separate treatment below.
\end{remark}
\begin{theorem}[Sonar-Radon relations in $\h{n}$] \label{t:MainResult}
The sonar transform $\sonar[f]$ determines $\radon[f]$ by means of
the following operator identities:
\begin{eqnarray}
    \radonh &=&
      \D{\frac{n-1}{2}} \circ \radonh \circ \sonar,
    \label{e:Horizontal} \\
    \radonv &=&
      \Lim  \circ \sonar,
    \label{e:Vertical} \\
    \radons &=&
      \W \circ \radonw[1/y] \circ
      \D{\frac{n-2}{2}} \circ \overline\radon \circ \sonar.
    \label{e:Slanted}
\end{eqnarray}
Here $\radonc$ stands for the centerset Radon transform; $\radonw$
stands for the weighted Radon transform \eqref{e:Radonw} from
Definition \ref{d:Radonw} in Section \ref{s:SlantedLines};
$\D{\nu}$ and $\W$ denote fractional operators defined in Section
\ref{s:FractionalCalculus} by \eqref{e:D} and \eqref{e:W},
respectively; $\Lim$  represents an infinite limit defined in
Section \ref{s:SlantedLines} by Equation \eqref{e:Lim}.
\end{theorem}
We organize the proof  as follows. Sections
\ref{s:FractionalCalculus} and \ref{s:PlaneWaves} contain
necessary analytical tools: Section \ref{s:FractionalCalculus}
details the fractional operators $\D{\nu}$ and $\W$ (and their
inverses); Section \ref{s:PlaneWaves} collates identities for
spherical integrals of plane waves from F.~John's classic
\textit{Plane Waves and Spherical Means} \cite{John-1981}, for use
in Section \ref{s:SlantedPlanes}. Section \ref{s:HorizontalPlanes}
and \ref{s:VerticalPlanes} treat $\radonh$ and $\radonv$,
respectively. Sections \ref{s:SlantedLines} and
\ref{s:SlantedPlanes} treat $\radons$ and thus complete the proof
of our sonar-Radon relations: Section \ref{s:SlantedLines}
motivates the choice of the weight $\weight=1/y$ for operator
$\radonw[1/y]$ and establishes results in dimension two; Section
\ref{s:SlantedPlanes} generalizes these results to higher
dimensions. Section \ref{s:Conclusion} analyzes the main result
and offers closing remarks.

\section{Fractional Calculus}
\label{s:FractionalCalculus}

We recall notions from fractional calculus, especially regarding
operators $\D{\nu}$ and $\W$ appearing in our sonar-Radon
relations.  The standing assumption that all fractional operators
act on the last variable of smooth compactly supported functions
will avoid those various delicate issues discussed at length by
Samko et al in \cite{Samko-et-al-1993}. Compact support  obviates
potential divergence; smoothness  ensures commutativity of
operators.
\par
As we choose to view the fractional integrals $\I{\nu}$  as the
fundamental fractional operators, we develop all other fractional
operators out of these.
\begin{definition} \label{d:I}
For
\begin{math}
  \nu > 0
\end{math}
we define $\I{\nu}$, a type of fractional integral, by
\begin{equation} \label{e:I}
  \I{\nu}[g](y) =
  \frac{2 \, \pi^{\nu}}{\Gamma(\nu)} \, y \,
  \int_{0}^{y}
    \left( y^{2} - s^{2} \right)^{\nu-1} \, g(s) \,
  ds.
\end{equation}
For convenience, we also set
\begin{math}
  \I{0} = \id.,
\end{math}
the identity operator.
\end{definition}
According to Lemma \ref{l:IComposition} below, the set
\begin{math}
  \left\{ \I{\nu} \right\}_{\nu=0}^{\infty}
\end{math}
of fractional integral operators forms a monoid (semigroup with
identity) under composition.
\begin{lemma} \label{l:IComposition}
The fractional integrals in Definition \ref{d:I} satisfy the
composition law
\begin{equation} \label{e:IComposition}
  \I{\mu} \circ \I{\nu} = \I{\mu + \nu},
\end{equation}
which holds for all
\begin{math}
  \mu, \nu \geq 0.
\end{math}
\end{lemma}
\begin{proof}
Equation \eqref{e:IComposition} certainly holds if either $\mu=0$
or $\nu=0$ on account of the convention $\I{0}=\id.$ So assume
$\mu,\nu>0$.
\par
Using Equation \eqref{e:I}, we express the composition
\begin{math}
  (\I{\mu} \circ \I{\nu})[g](y)
\end{math}
as an iterated integral
\begin{equation} \label{e:IComposition1}
  \frac{2 \, \pi^{\mu}}{\Gamma(\mu)} \,
  \frac{2 \, \pi^{\nu}}{\Gamma(\nu)} \, y \,
  \int_{0}^{y}
    \left( y^{2} - s^{2} \right)^{\mu-1} \, s \,
    \left[
      \int_{0}^{s}
        \left( s^{2} - t^{2} \right)^{\nu-1} \, g(t) \,
      dt
    \right] \,
  ds.
\end{equation}
Changing the order of integration in \eqref{e:IComposition1}
yields the following expression for
\begin{math}
  (\I{\mu} \circ \I{\nu})[g](y):
\end{math}
\begin{equation} \label{e:IComposition2}
  \frac{2 \, \pi^{\mu + \nu}}
       {\Gamma(\mu) \, \Gamma(\nu)} \,
  y \,
  \int_{0}^{y}
    \left[  I(t)
    \right] \,
    g(t) \,
  dt.
\end{equation}
where
\begin{equation} \label{e:IComposition2InsideInt1}
I(t)=  \int_{t}^{y}
    \left( y^{2} - s^{2} \right)^{\mu-1} \,
    \left( s^{2} - t^{2} \right)^{\nu-1} \, 2s \,
  ds.
\end{equation}
Substituting
\begin{math}
  s^{2} = \left( y^{2} - t^{2} \right) \, p + t^{2},
\end{math}
gives
\begin{displaymath}
  I(t) =
 C \left( y^{2} - t^{2} \right)^{\mu + \nu -1}
\end{displaymath}
with the constant
\begin{displaymath}
  C=\int_{0}^{1}
    (1 - p)^{\mu-1} \, p^{\nu-1} \,
  dp
\end{displaymath}
taking the form of Euler's integral of the first kind (see
\cite{Erdelyi-et-al-1981} p.948 \# 8.380.1) with value given by
\begin{displaymath}
  C =
  \frac{\Gamma(\mu) \, \Gamma(\nu)}
       {\Gamma(\mu + \nu)}.
\end{displaymath}
Thus
\begin{equation} \label{e:IComposition3}
(\I{\mu} \circ \I{\nu})[g](y)=
  \frac{2 \, \pi^{\mu + \nu}}
       {\Gamma(\mu) \, \Gamma(\nu)} \,
  y \,
  \int_{0}^{y}
    \left[
    \frac{\Gamma(\mu) \, \Gamma(\nu)}{\Gamma(\mu + \nu)}
    \left( y^{2} - t^{2} \right)^{\mu + \nu - 1} \,
    \right] \,
    g(t) \,
  dt.
\end{equation}
Cancelling four gamma terms, we obtain
\begin{displaymath}
  (\I{\mu} \circ \I{\nu})[g](y)=\frac{2 \, \pi^{\mu + \nu }}
       {\Gamma(\mu + \nu )} \,
  y \,
  \int_{0}^{y}
    \left( y^{2} - t^{2} \right)^{\mu + \nu - 1} \, f(t) = \I{\mu + \nu}[g](y)\,
  dt,
\end{displaymath}
as desired.
\end{proof}
The justification for the terminology \textit{fractional
integrals} for the $\I{\nu}$ rests on the semigroup property and
the observation that
\begin{displaymath}
  \I{1}[g](y) = 2 \, \pi \, y
  \int_{0}^{y}
    g(s) \,
  ds,
\end{displaymath}
a \textit{scaled} antiderivative. Now $\I{1}$ admits a left
inverse in the form
\begin{equation} \label{e:D1}
  \D{1}[g](y) := \frac{1}{2 \, \pi} \,
  \frac{d}{dy}
  \left[
    \frac{g(y)}{y}
  \right].
\end{equation}
So, symbolically, we have
\begin{math}
  \D{1} \circ \I{1} = \I{0} = \id,
\end{math}
and, more generally,
\begin{lemma} \label{l:D0}
For all $\nu>1$,
\begin{math}
  \D{1} \circ \I{\nu} = \I{\nu-1}.
\end{math}
\end{lemma}
\begin{proof}
Using Lemma \ref{l:IComposition} and the associativity of operator
composition
\begin{displaymath}
  \D{1} \circ \I{\nu} =
  \D{1} \circ (\I{1} \circ \I{\nu-1}) =
  (\D{1} \circ \I{1}) \circ \I{\nu-1} =
  \I{\nu-1}.
\end{displaymath}
\end{proof}
We will now use Lemmas \ref{l:IComposition} and \ref{l:D0} to
construct fractional derivatives of arbitrary order.
\begin{definition} \label{d:D}
For
\begin{math}
  \nu > 0
\end{math}
the fractional derivative $\D{\nu}$ is defined in terms of
\eqref{e:I} and \eqref{e:D1} as a mapping
\begin{equation} \label{e:D}
\D{\nu} =
   \D{1}^{\ceil{\nu}} \circ \I{\ceil{\nu}-\nu},
\end{equation}
where $\ceil{\nu}$ is the smallest integer greater than or equal
$\nu$.  For consistency with Definition \ref{d:I}, we define
\begin{math}
  \D{0} = \id.
\end{math}
\end{definition}
As an immediate consequence of Lemmas \ref{l:IComposition} and
\ref{l:D0}, we have the following corollary.
\begin{corollary} \label{c:iI}
For
\begin{math}
  \nu \geq 0
\end{math}
the fractional derivative $\D{\nu}$ is the inverse of the
fractional integral $\I{\nu}$, i.e.:
\begin{math}
  \D{\nu} \circ \I{\nu} = \id.
\end{math}
\end{corollary}
In Section \ref{s:SlantedLines}, we will encounter a fractional
operator $\V$ and require its inverse $\W$ to deduce Equation
\eqref{e:Slanted} in Theorem \ref{t:MainResult}.
\begin{definition} \label{d:V}
Set
\begin{equation} \label{e:V}
  \V[g](\beta):=
  \int_{0}^{\beta}
    \frac{2 \, \sin \beta}{\sqrt{\sin^{2} \beta - \sin^{2} \theta}} \,
    g(\theta) \,
 d\theta.
\end{equation}
\end{definition}
From the definitions of $\I{\frac{1}{2}}$ and $\V$
\begin{displaymath}
  \I{\frac{1}{2}}[g](\sin \beta) =
  \V[\cos\beta \, g(\sin \beta)](\beta).
\end{displaymath}
We now cast this identity of functions as an identity of
operators.  Define
\begin{math}
  {\cal Q}[g](\beta) := g(\sin\beta)
\end{math}
and
\begin{math}
  {\cal K}[g](\beta) := \cos\beta \cdot g(\sin\beta).
\end{math}
Then the identity above says
\begin{displaymath}
  {\cal Q} \circ \I{\frac{1}{2}} = \V \circ {\cal K}.
\end{displaymath}
Thus
\begin{displaymath}
  \W := \V^{-1} =
  {\cal  K} \circ \D{\frac{1}{2}} \circ {\cal Q}^{-1}.
\end{displaymath}
\begin{lemma} \label{l:W}
\begin{equation} \label{e:W}
   \W[g](\beta):=
    \frac{1}{\pi} \,
    \frac{d}{d \beta}
    \int_{0}^{\beta}
      \frac{\cos \theta}{\sqrt{\sin^{2} \beta - \sin^{2} \theta}} \,
      g(\theta) \,
    d\theta
\end{equation}
\end{lemma}
\begin{proof}
We apply the composition of operators
\begin{math}
  {\cal  K} \circ \D{\frac{1}{2}} \circ {\cal Q}^{-1}
\end{math}
to a function $g$.  First,
\begin{displaymath}
  {\cal Q}^{-1}[g](s) = g(\arcsin s).
\end{displaymath}
Next we use Equation \eqref{e:D} to write
\begin{equation*}
  \begin{split}
    ({\D{\frac{1}{2}}\circ \cal Q}^{-1})[g](\beta) &=
    (\D{1} \circ \I{\frac{1}{2}} \circ {\cal Q}^{-1})[g](\beta)\\
    &=\frac{1}{2 \, \pi} \,
    \frac{d}{d\beta}
    \left[
      \frac{1}{\beta } \cdot 2 \, \beta
      \int_{0}^{\beta }
        \frac{g(\arcsin s)}{\sqrt{\beta^{2}-s^{2}}} \,
      ds
    \right]\\
    &=\frac{1}{\pi} \,
    \frac{d}{d\beta}
    \int_{0}^{\beta}
      \frac{g(\arcsin s)}{\sqrt{\beta^{2}-s^{2}}} \,
    ds
  \end{split}
\end{equation*}
Finally, we apply the operator $K$, which replaces $\beta$ with
$\sin\beta$ and multiplies the result by $\cos\beta$, to get
\begin{equation*}
  \begin{split}
    ({K \circ \D{\frac{1}{2}}\circ \cal Q}^{-1})[g](\beta)
    &=
    \cos\beta \, \frac{1}{\pi} \,
    \frac{d}{d(\sin\beta)}
    \int_{0}^{\sin \beta}
      \frac{g(\arcsin s)}{\sqrt{\sin^{2}\beta-s^2}} \,
    ds
  \end{split}
\end{equation*}
whereupon the substitution $s = \sin\theta$ yields the statement
of the lemma.
\end{proof}

\section{Plane waves and spherical means}
\label{s:PlaneWaves}

The identities for integrals over spheres and balls in $\r{n}$
collected here, combined with the formulae from Section
\ref{s:FractionalCalculus}, form the crux of the derivations
presented in Sections \ref{s:HorizontalPlanes} and
\ref{s:SlantedPlanes}. In particular, we state the Co-area
Formula, following \cite{Evans-Gariepy-1992}, and develop some of
its consequences.
\begin{theorem}[Co-area formula] \label{t:coarea}
Let
\begin{math}
  u: \r{n} \mapsto \r{}
\end{math}
be Lipschitz continuous and assume that for almost every
\begin{math}
  r \in \r{}
\end{math}
the level set
\begin{displaymath}
  \left\{ x \in \r{n} \, | \, u(x) = r \right\}
\end{displaymath}
is a smooth, $(n-1)$-dimensional surface in $\r{n}$.  Suppose that
also
\begin{math}
  g: \r{n} \mapsto \r{}
\end{math}
is continuous and locally integrable.   Then
\begin{equation} \label{e:coarea}
  \int_{x \in \r{n}}
    g \cdot |\grad u| \,
  dx =
  \int_{-\infty}^{\infty}
  \left[
    \int_{\{ u = r \}}
      g \,
    dS
  \right] \,
  dr,
\end{equation}
where $dS$ denotes surface measure on the level set
\begin{math}
  \{ u = r \}.
\end{math}
\end{theorem}
By setting
\begin{math}
  u = |x|
\end{math}
in Theorem \ref{t:coarea}, one obtains a standard identity for
converting integrals over balls into integrals over spheres in
$\r{n}$ which we state in Lemma \ref{l:Polar}.
\begin{lemma}[Polar Coordinates] \label{l:Polar}
Let
\begin{math}
  g: \b{r}{n} \mapsto \r{}
\end{math}
be a continuous function on a ball of radius $r$ in Euclidean
space. Then
\begin{equation} \label{e:Polar}
  \int_{|x| < r}
    g(x) \,
  dx =
  \int_{0}^{r}
  \left[
    \int_{x \in \s{p}{n-1}}
      f(x) \,
    \dS{x}
  \right] \,
  dp,
\end{equation}
where $\dS{x}$ denotes surface measure on the sphere $\s{p}{n-1}$
of radius $p$.
\end{lemma}
Differentiating Equation \eqref{e:Polar}  gives rise to the
following.
\begin{corollary} \label{c:Polar1}
Let
\begin{math}
  g: \b{R}{n} \mapsto \r{}
\end{math}
be a continuous function. Then
\begin{equation} \label{e:Polar1}
  \frac{d}{dr} \,
  \int_{x \in \b{r}{n}}
    g(x) \,
  dx =
  \int_{x \in \s{r}{n-1}}
    g(x) \,
  \dS{x}
\end{equation}
holds for all
\begin{math}
  0 < r < R.
\end{math}
\end{corollary}
Often it is convenient to replace integration over a sphere of
radius $r$ with integration over a unit sphere.  As we remark
below, this can be accomplished through a simple substitution.
\begin{remark} \label{r:Polar2}
Let
\begin{math}
  g: \s{r}{n} \mapsto \r{}
\end{math}
be a continuous function on a sphere of radius $r$.  Then
\begin{equation} \label{e:Polar2}
  \int_{x \in \s{r}{n-1}}
    g(x) \,
  \dS{x} =
  r^{n - 1} \,
  \int_{\theta \in \s{1}{n-1}}
    g(r \, \theta) \,
  \ds{\theta},
\end{equation}
where $\ds{\theta}$ is the surface measure on a unit sphere.
\end{remark}
Using Remark \ref{r:Polar2}, we  recast Equation \eqref{e:Polar}
in the form we shall find most useful:
\begin{equation} \label{e:polar}
  \int_{|x| < r}
    g(x) \,
  dx =
  \int_{0}^{r}
  s^{n - 1} \,
  \left[
    \int_{\theta \in \s{1}{n-1}}
      g(s \, \theta) \,
    \ds{\theta}
  \right] \,
  ds.
\end{equation}
\par
Throughout the rest of this section $g$ denotes a continuous
function of a scalar variable. Fix
\begin{math}
  v \in \r{n}.
\end{math}
Following F.~John in \cite{John-1981}, we call $G(x) := g(v \cdot
x)$ a \textit{plane wave} with normal $v$; such a function is
constant on planes perpendicular to $v$.
\par
We follow F.~John in \cite{John-1981}  to reduce integrals of
plane waves over spheres and balls  to single-dimensional
integrals.  On the plane
\begin{math}
  v \cdot x = |v| \, p
\end{math}
the plane wave $g(v \cdot x)$ has constant value $g(|v| \, p)$;
the  intersection of that plane with the ball of radius $r$ forms
a ball of radius $\sqrt{r^2-p^2}$. Thus
\begin{equation} \label{e:John1.1}
  \begin{split}
    \int_{|x| < r}
      g(v \cdot x) \,
    dx
    &=  \,
    \int_{-r}^{+r} \makebox{\rm Vol}(\b{\sqrt{r^2-p^2}}{n-1})
      g(|v| \, p) \,
    dp\\
    &=
    \frac{|\s{1}{n-2}|}{n-1} \,
    \int_{-r}^{+r}
      \left( r^{2} - p^{2} \right)^{\frac{n - 1}{2}} \,
      g(|v| \, p) \,
    dp.
  \end{split}
\end{equation}
where $|\s{1}{n-2}|$ denotes the total surface measure of a unit
sphere in $\r{n-1}$.    Differentiation of Equation
\eqref{e:John1.1} with respect to $r$ followed by evaluation at
\begin{math}
  r = 1,
\end{math}
leads to the following fundamental identity (c.f.
\cite{John-1981}, p.8).
\begin{theorem} \label{t:John}
Let
\begin{math}
  v \in \r{n}
\end{math}
be a fixed vector and let
\begin{math}
  g: \mathbb{R} \mapsto \mathbb{R}
\end{math}
be a continuous function.  Then
\begin{equation} \label{e:John}
  \int_{\theta \in \s{1}{n-1}}
    g(v \cdot \theta) \,
  \ds{\theta} =
  |\s{1}{n-2}| \,
  \int_{-1}^{+1}
    \left( 1 - p^{2} \right)^{\frac{n - 3}{2}} \,
    g(|v| \, p) \,
  dp.
\end{equation}
\end{theorem}
We state two consequences of Theorem \ref{t:John}. If $v$ has unit
length, then Equation \eqref{e:John} becomes an identity for
``spherical plane waves'' :
\begin{corollary} \label{c:John1}
Let
\begin{math}
 \omega \in \s{1}{n-1}
\end{math}
be a fixed point on a unit sphere in $\r{n}$ and let
\begin{math}
  g: \mathbb{R} \mapsto \mathbb{R}
\end{math}
be a continuous function.  Then
\begin{equation} \label{e:John1}
  \int_{\theta \in \s{1}{n-1}}
    g(\omega \cdot \theta) \,
  \ds{\theta} =
  |\s{1}{n-2}| \,
  \int_{-1}^{+1}
    \left( 1 - p^{2} \right)^{\frac{n-3}{2}} \,
    g(p) \,
  dp.
\end{equation}
\end{corollary}
Alternatively, on setting $g = 1$, Theorem \ref{t:John} gives a
recursion for  the total measure of a unit sphere (c.f.
\cite{John-1981},  p.9). From this recursion follows the
well-known surface area formula:
\begin{equation} \label{e:John2}
  |\s{1}{n-1}| =
  \frac{2 \, \pi^{\frac{n}{2}}}
       {\Gamma \left( \frac{n}{2} \right)}.
\end{equation}
We shall use this formula to connect our fractional integrals with
geometric transforms.

\section{Integrals over horizontal planes}
\label{s:HorizontalPlanes}

We now prove Equation \eqref{e:Radonh} of our Main Theorem:
\begin{displaymath}
  \radonh = \D{\frac{n-1}{2}} \circ \radonh \circ \sonar.
\end{displaymath}
Since $\D{\nu}$ inverts the fractional integral $\I{\nu}$,  it
suffices to show that
\begin{equation} \label{e:HP}
  \radonh \circ \sonar =
  \I{\frac{n - 1}{2}} \circ \radonh.
\end{equation}
We compare
\begin{math}
  (\radonh \circ \sonar)[f]
\end{math}
with
\begin{math}(\I{\frac{n - 1}{2}} \circ \radonh)[f]\end{math} for
\begin{math}
  f \in \cc{\h{n}}.
\end{math}
Directly from the definitions of the operators $\sonar$ and
$\radonh$ (Equations \eqref{e:sonar},\eqref{e:Radonh})
\begin{displaymath}
(\radonh \circ \sonar)[f](y) =
  \int_{x \in \r{n-1}}
  \left[
    \int_{|t| < y}
      f \left( x + t, \sqrt{y^{2} - |t|^{2}} \right) \,
    \frac{y \, dt}{\sqrt{y^{2} - |t|^{2}}}
  \right] \,
  dx
\end{displaymath}
which, upon interchanging the order of integration, equals
\begin{displaymath}
  \int_{|t| < y}
  \left[
    \int_{x \in \r{n-1}}
      f \left( x + t, \sqrt{y^{2} - |t|^{2}} \right) \,
    dx
  \right] \,
  \frac{y \, dt}{\sqrt{y^{2} - |t|^{2}}}.
\end{displaymath}
We recognize the inner integral as $\radonh[f](\sqrt{y^{2} -
|t|^{2}})$ and deduce
\begin{equation} \label{e:HP3}
  \left( \radonh \circ \sonar \right)[f](y) =
  \int_{|t| < y}
    \radonh[f] \left( \sqrt{y^{2} - |t|^{2}} \right) \,
  \frac{y \, dt}{\sqrt{y^{2} - |t|^{2}}}.
\end{equation}
The right-hand side of Equation \eqref{e:HP3} is an integral over
a ball in $\r{n-1}$ whose integrand
\begin{displaymath}
  \radonh[f] \left( \sqrt{y^{2} - |t|^{2}} \right) \,
  \frac{y}{\sqrt{y^{2} - |t|^{2}}}
\end{displaymath}
is a radial function of $t$.  Therefore, using Equation
\eqref{e:polar} from Section \ref{s:PlaneWaves}, we  rewrite
Equation \eqref{e:HP3} in the form
\begin{displaymath}
  |\s{1}{n-2}| \, y \,
  \int_{0}^{y}
      \frac{\radonh[f] \left( \sqrt{y^{2} - r^{2}} \right)}
           {\sqrt{y^{2} - r^{2}}} \,
       r^{n-2} \,
  dr,
\end{displaymath}
which, after substituting
\begin{math}
  r = \sqrt{y^{2} - s^{2}}
\end{math}
becomes
\begin{equation} \label{e:HP5}
   |\s{1}{n-2}| \, y \,
  \int_{0}^{y}
    \left( y^{2} - s^{2} \right)^{\frac{n - 3}{2}}
    \radonh[f](s)\,ds =(\I{\frac{n - 1}{2}} \circ \radonh)[f]
  \end{equation}
by Definition \ref{d:I} in Section \ref{s:FractionalCalculus} and
\eqref{e:John2}.

\section{Integrals over vertical planes}
\label{s:VerticalPlanes}

\begin{figure}[ht!]
  \begin{center}
    \scalebox{.35}{\includegraphics{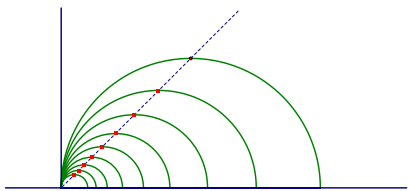}} \\
    \caption{Vertical rays as limits of arcs}\label{fig:vertical}
  \end{center}
\end{figure}

Figure \ref{fig:vertical} suggests viewing a vertical hyperplane
as a limiting case of expanding tangent spheres with a fixed point
of tangency located on the centerset.

Accordingly,
\begin{displaymath}
  \radonv[g](\omega,p) =
  \lim_{|s|\to\infty}
  \sonar[g](\omega \, s, |s-p|).
\end{displaymath}
If we make the definition
\begin{equation} \label{e:Lim}
  \Lim [f](\omega,p) :=
  \lim_{|s| \to \infty} f(\omega  s, |p-s|).
\end{equation}
then
\begin{math}
  \radonv = \Lim  \circ \sonar,
\end{math}
as desired.

\section{Integrals over slanted lines}
\label{s:SlantedLines}

In $\h{2}$, our desired sonar-Radon relation Equation
\eqref{e:Slanted} reduces to:
\begin{equation} \label{e:SR}
  \radons = \W \circ \radonw[1/y] \circ \radonc \circ \sonar.
\end{equation}
(As
\begin{math}
  \D{\frac{n - 2}{2}} = \D{0} = \id,
\end{math}
no fractional derivative appears.) Below, this formula emerges as
the  foundation for the general case.
\par

\par
In dimension two, the centerset has dimension one.  As hyperplanes
in dimension one coincide with points, just in this section we
will encode them as such (rather than as pairs of a unit vector
and a magnitude). With this encoding the centerset Radon transform
$\radonc$ reduces to the identity map, so we must prove that
\begin{equation} \label{e:SR1}
  \radons = \W \circ \radonw[1/y]  \circ \sonar.
\end{equation}
Applying $\V$, the inverse of $\W$ to both sides yields the
equivalent statement
\begin{equation} \label{e:SR2}
  \V \circ  \radons = \radonw[1/y] \circ \sonar
\end{equation}
for which we will aim.
\par
Dimension two affords us a simple formula for the sonar transform:
\begin{equation} \label{e:S2}
    \sonar[f](x,y) =
    \int_{0}^{\pi}
      f(x + y \, \cos \phi, y \, \sin \phi) \, y \,
    d\phi.
\end{equation}
We now furnish the definition of $\radonw[\sigma]$ --- a type of
weighted  Radon transform on $\h{2}$.
\begin{definition} \label{d:Radonw}
Fix any set $T$.  Consider a function $g$ on $T \times \h{2}$ such
that
\begin{math}
  g(\omega,\cdot) \in \cc{\h{2}}
\end{math}
for each $\omega$ in $T$. For
\begin{math}
  \beta \in (0, \pi / 2)
\end{math}
and non-negative weight
\begin{math}
  \weight: \h{} \mapsto \h{},
\end{math}
define a weighted Radon transform  by
\begin{equation} \label{e:Radonw}
  \radonw[\sigma][g](\omega,p,\beta) =
  \int_{0}^{\infty}
    g \left(\omega, p + s, s \, \sin \beta \right) \,
    \weight(s) \,
  ds.
\end{equation}
\end{definition}
This section has $T$ a singleton and we thus suppress the variable
$\omega$.
\par
\begin{figure} [ht!]
  \begin{center}
    \scalebox{.35}{\includegraphics{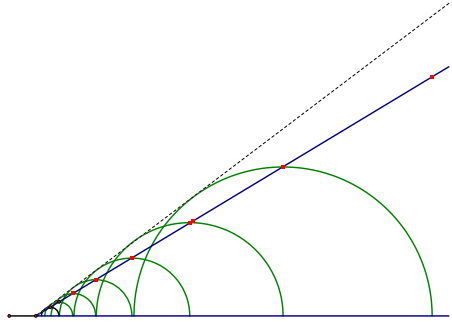}} \\
    \caption{Semicircles tangential to a ray} \label{fig:slanted}
  \end{center}
\end{figure}

Consider the composition
\begin{math}
  \radonw[\sigma] \circ \sonar
\end{math}
for a general weight $\weight$. $\sonar[f](p,y)$ means the
integral of $f$ over over a radius $y$ circle centered at $p$ on
the $x$-axis.  By definition, the operator $\radonw[\sigma]$
integrates functions along rays with slope $\sin\beta$.  This
makes
\begin{math}
  \radonw[\weight][\sonar[f]\,](p,\beta)
\end{math}
a weighted integral of integrals of $f$ over a family of circles,
as in Figure \ref{fig:slanted}.
\par
One also sees from the figure that arcs with apexes on slanted
rays sweep infinite wedges. If the apexes lie on a ray with slope
$\sin\beta$ then the corresponding wedge has slope $\tan\beta$.
Therefore
\begin{math}
  \radonw[\weight][\sonar[f]\,](p,\beta)
\end{math}
can be expressed as an integral over an infinite wedge with vertex
at $p$ on the $x$-axis and angular measure $\beta$.  Explicitly,
by means of Equations \eqref{e:S2} and \eqref{e:Radonw},
\begin{equation} \label{e:AS}
  \begin{split}
    &\radonw[\weight][\sonar[f]\,](p,\beta)\\&=
  \int_{0}^{\infty}
  \left[
    \int_{0}^{\pi}
      f
      \left(
        p + t \, (1 + \sin \beta \, \cos \phi) ,
        t \, \sin \beta \, \sin \phi
      \right) \,
    d\phi
  \right] \,
  t \, \sin \beta \, \weight(t) \,
  dt
  \end{split}
\end{equation}
We now make a change of variables designed to simplify the
argument of $f$ in Equation \eqref{e:AS}. Define
\begin{math}
  \Psi: \r{2} \rightarrow \r{2}
\end{math}
 by
\begin{displaymath}
  \left[
    \begin{array}{c}
      \phi \\
         t \\
    \end{array}
  \right] \stackrel{\Psi}{\mapsto}
  \left[
    \begin{array}{c}
      p + t \, (1 + \sin \beta \, \cos \phi) \\
                     t \, \sin \beta \, \sin \phi \\
    \end{array}
  \right].
\end{displaymath}
Observe that $\Psi$ sends a line segment connecting $(0,t)$ and
$(\pi,t)$ to the semicircle centered at $(p+t,0)$ with radius
$t\sin\beta$.
\par
On each of the two semi-infinite strips
\begin{displaymath}
  \left( 0, \frac{\pi}{2} + \beta \right) \times (0,\infty)
  \quad \text{and} \quad
  \left( \frac{\pi}{2} + \beta, \pi \right) \times (0,\infty)
\end{displaymath}
$\Psi$ acts as a diffeomorphism to the infinite wedge
\begin{displaymath}
  \left\{
    (x,y) \in \r{2} \mid
    p < x,
    0 < y < \tan \beta \, (x - p)
  \right\}
\end{displaymath}
shown on Figure \ref{fig:slanted}.  In terms of $\Psi$, the double
integral \eqref{e:AS} over a wedge can be written as a sum of two
integrals over infinite strips
\begin{equation} \label{e:ASPsi}
  \begin{split}
    &\int_{0}^{\infty}
      \left[
        \int_{0}^{\frac{\pi}{2} + \beta}
          (f \circ \Psi)(\phi,t) \,
        d\phi
      \right] \, t \, \sin \beta \, \weight(t) \,
    dt \\&+
    \int_{0}^{\infty}
      \left[
        \int_{\frac{\pi}{2} + \beta}^{\pi}
          (f \circ \Psi)(\phi,t) \,
        d\phi
      \right] \, t \, \sin \beta \, \weight(t) \,
    dt.
  \end{split}
\end{equation}
Introducing polar coordinates $(\rho,\theta)$,
\begin{math}
  \rho \in (0, \infty)
\end{math},
\begin{math}
  \theta \in (0, \beta)
\end{math}
in the wedge
\begin{displaymath}
  (p + t \, (1 + \sin\beta \, \cos\phi) ,  t \, \sin\beta \, \sin\phi) =
  (p + \rho \, \cos\theta, \rho \, \sin\theta)
\end{displaymath}
gives us the relations:
\begin{align}
  t \, (1 + \sin \beta \, \cos \phi) &=
     \rho \, \cos \theta \label{e:r2r1polarx}\\
  t \, \sin \beta \, \sin \phi &=
    \rho \, \sin \theta. \label{e:r2r1polary}
\end{align}
Using \eqref{e:r2r1polarx} and \eqref{e:r2r1polary}, we shall now
change \eqref{e:ASPsi} into a much more amenable expression.
\par
In order to transform \eqref{e:ASPsi}, we need to express the old
variables $(\phi,t)$ in terms of the new variables $(\theta,\rho)$
and find the corresponding Jacobian factors: one for each integral
in \eqref{e:ASPsi}.  From the algebraic point of view, it is
easier to find $\phi$. Divide Equation \eqref{e:r2r1polary} by
Equation \eqref{e:r2r1polarx}: this eliminates variables $\rho$
and $t$. Next use trigonometric identities to solve the resulting
relation between angles as follows:
\begin{equation} \label{e:r2r1phi}
  \begin{split}
    \frac{\operatorname{Eq.} \eqref{e:r2r1polary}}
         {\operatorname{Eq.} \eqref{e:r2r1polarx}}
   &\Rightarrow
      \frac{\sin \beta \, \sin \phi}
           {1 + \sin \beta \, \cos \phi} =
      \tan \theta
    \\&
    \Rightarrow
      \sin \beta \, \sin \phi \, \cos \theta =
        \sin \theta +  \sin \beta \, \cos \phi \, \sin \theta
    \\&
    \Rightarrow
      \sin \phi \, \cos \theta - \cos \phi \, \sin \theta =
        \frac{\sin \theta}{\sin \beta}
    \\&
    \Rightarrow
      \sin (\phi - \theta) =
        \frac{\sin \theta}{\sin \beta}
    \\&
    \text{we get two solutions:} \\&
    \left\{
      \begin{split}
        \phi_{1} &= \theta + \sin^{-1}
        \left(
          \frac{\sin \theta}{\sin \beta}
        \right), \\
        \phi_{2} &= \theta + \pi - \sin^{-1}
        \left(
          \frac{\sin \theta}{\sin \beta}
        \right)
      \end{split}
    \right\}.
  \end{split}
\end{equation}
With these expressions for the angular variable $\phi$, we may now
find the corresponding values of $t$ as outlined  in Equation
\eqref{e:r2r1t} below:
\begin{equation} \label{e:r2r1t}
  \begin{split}
    \operatorname{Eq.} \eqref{e:r2r1polarx} \times \cos \theta &+
    \operatorname{Eq.} \eqref{e:r2r1polary} \times \sin \theta \\&
    \Rightarrow
      t \, \cos \theta + t \, \sin \beta \,
        (\cos \phi \, \cos \theta + \sin \phi \, \sin \theta) = \rho
    \\&
    \Rightarrow
      t \, \cos \theta + t \, \sin \beta \, \cos (\phi - \theta) = \rho
    \\&
    \Rightarrow
      t =
        \frac{\rho}{\cos \theta + \sin \beta \, \cos (\phi - \theta)}
    \\&
    \text{using the identity following from Equation \eqref{e:r2r1phi}:}\\&
    \left(
      \begin{split}
            \sin \beta \, \cos (\phi - \theta) =
        \pm \sin \beta \, \sqrt{1 - \sin^{2} (\phi - \theta)} &=\\
        \pm \sin \beta \, \sqrt{1 - (\sin^{2} \theta / \sin^{2} \beta)} =
        \pm \sqrt{\sin^{2} \beta - \sin^{2} \theta}
      \end{split}
    \right) \\&
    \text{we get two solutions:} \\&
    \left\{
      \begin{split}
        t_{1} &=
          \frac{\rho}
               {\cos \theta + \sqrt{\sin^{2} \beta - \sin^{2}
               \theta}}, \\
        t_{2} &=
          \frac{\rho}
               {\cos \theta - \sqrt{\sin^{2} \beta - \sin^{2} \theta}}
      \end{split}
    \right\}.
  \end{split}
\end{equation}
\begin{remark} \label{r:geometricJ}
The relation between  $(\phi,t)$ and $(\theta,\rho)$ can also be
derived geometrically.
\par
\begin{figure} [ht!]
  \begin{center}
    \scalebox{.75}{\includegraphics{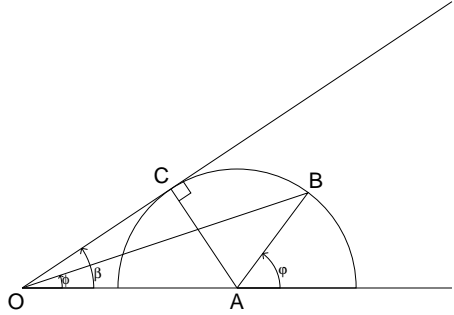}} \\
    \caption{Polar coordinates} \label{fig:polar}
  \end{center}
\end{figure}

Consider semicircles inscribed in a fixed wedge of angular measure
$\beta$ as in Figure \ref{fig:polar}. Our old variables $(\phi,t)$
specify a semicircle and then a point $B$ on it: from $t=OA$ we
learn the center of the semicircle, and then from the tangency
also its radius;  $\phi$ locates the point $B$ since $\angle OAB =
\pi - \phi$. From the right triangle $\triangle OAC$, we find the
radius of the semicircle
\begin{math}
  r = t \, \sin\beta.
\end{math}
A ray issuing from $O$ at angle $\theta<\beta$ with $OA$ will meet
the semicircle twice and we take $B$ as the second intersection.
The Law of Sines applied to
 triangle $\triangle OAB$ gives:
\begin{displaymath}
  \frac{\sin (\angle AOB)}{AB} =
  \frac{\sin (\angle OBA)}{OA}
  \Rightarrow
  \frac{\sin\theta}{t \, \sin\beta} =
  \frac{\sin(\phi - \theta)}{t},
\end{displaymath}
which is equivalent to \eqref{e:r2r1phi}.  Then the Law of
Cosines, in the form
\begin{displaymath}
  OA^{2} + AB^{2} - 2 \, OA \, AB \, \cos(\angle OAB) = OB^{2},
\end{displaymath}
tells us
\begin{displaymath}
  t^{2} + t^{2} \, \sin^{2}\beta - 2 \, t^{2} \, \sin\beta \,
  \cos(\pi - \phi) = \rho^{2}.
\end{displaymath}
After simplification \eqref{e:r2r1t} results.
\end{remark}
As follows from Equation \eqref{e:r2r1phi} the angular variable
$\phi$ does not depend on $\rho$.  Therefore the $2$-by-$2$
Jacobian matrix is triangular and its determinant is given by
\begin{displaymath}
  \det
  \left[
    \begin{array}{cc}
      \frac{\partial \phi}{\partial \theta} & 0 \\
      \frac{\partial t}{\partial \theta} & \frac{\partial t}{\partial \rho} \\
    \end{array}
  \right] =
    \frac{\partial \phi}{\partial \theta} \times
    \frac{\partial t}{\partial \rho}.
\end{displaymath}
The values of the partial derivatives
\begin{math}
  \frac{\partial \phi_{i}}{\partial \theta}
\end{math}
and
\begin{math}
  \frac{\partial t_{i}}{\partial \rho}
\end{math}
for $i=1,2$ can be found through straightforward differentiation:
\begin{align*}
  \frac{\partial \phi_{i}}{\partial \theta} &=
    1 + \frac{(-1)^{i}}{\sqrt{1 - (\sin^{2} \theta / \sin^{2}
    \beta)}} \, \frac{\cos \theta}{\sin \beta} =
    \frac{\sqrt{\sin^{2} \beta - \sin^{2} \theta} + (-1)^{i} \cos \theta}
         {\sqrt{\sin^{2} \beta - \sin^{2} \theta}}, \\
  \frac{\partial t_{i}}{\partial \rho} &=
    \frac{1}{\cos \theta + (-1)^{i} \sqrt{\sin^{2} \beta - \sin^{2}
    \theta}}, \quad i = 1,2,
\end{align*}
whence follows that for both sets of variables $(\phi_{1},t_{i})$
$i=1,2$ the absolute value of the determinant of the Jacobian is
given by the same simple expression
\begin{displaymath}
  \frac{1}{\sqrt{\sin^{2} \beta - \sin^{2} \theta}}.
\end{displaymath}
We conclude that \eqref{e:ASPsi} can be written as a single
integral of the form
\begin{equation} \label{e:r2r1G}
  \int_{0}^{\beta}
    \left[
      \int_{0}^{\infty}
        \frac{ \sin \beta \,
        \left\{
          t_{1} \, \weight(t_{1}) + t_{2} \, \weight(t_{2})
        \right\}}
             {\sqrt{\sin^{2} \beta - \sin^{2} \theta}} \,
        f(p + \rho \, \cos \theta, \rho \, \sin \theta) \,
      d\rho
    \right] \,
  d\theta,
\end{equation}
where the values of $t_{i}$, $i = 1,2$ in the numerator are given
by \eqref{e:r2r1t}.
\par
Integral \eqref{e:r2r1G}, representing the composition
\begin{math}
  \left( \radonw \circ \sonar \right)
\end{math}
applied to $f$, becomes particularly simple if one sets the weight
\begin{math}
  \weight = 1 / y:
\end{math}
\begin{equation} \label{e:r2r1Gint}
  \begin{split}
    &\int_{0}^{\beta}
      \left[
        \int_{0}^{\infty}
          \frac{2 \, \sin \beta }
               {\sqrt{\sin^{2} \beta - \sin^{2} \theta}} \,
          f(p + \rho \, \cos \theta, \rho \, \sin \theta) \,
        d\rho
      \right] \,
    d\theta \\
    &=
    \int_{0}^{\beta}
      \frac{2 \, \sin \beta }
           {\sqrt{\sin^{2} \beta - \sin^{2} \theta}} \,
        \left[
          \int_{0}^{\infty}
            f(p + \rho \, \cos \theta, \rho \, \sin \theta) \,
          d\rho
        \right] \,
    d\theta.
  \end{split}
\end{equation}
Recognizing the bracketed integral as $\radons$ (via Equation
\eqref{e:Radons} from Section \ref{s:MainResult}) and noticing
that the outside integral is the fractional operator $\V$
(Definition \ref{d:V} from Section \ref{s:FractionalCalculus}), we
get
\begin{equation*}
  \begin{split}
    (\radonw[1/y] \circ \sonar)[f](p,\beta)&=
    \int_{0}^{\beta}
      \frac{2 \, \sin \beta }
           {\sqrt{\sin^{2} \beta - \sin^{2} \theta}} \,
           \radons{f}(p,\theta) \,
    d\theta\\
    &=
 (\V \circ \radons)[f](p,\beta),
  \end{split}
\end{equation*}
as desired.

\section{Integrals over slanted planes}
\label{s:SlantedPlanes}
%
We now prove in all dimensions the sonar-Radon relation
\eqref{e:Slanted} first stated in Theorem \ref{t:MainResult} and
reproduced below:
\begin{equation*}
  \radons = \W \circ \radonw[1/y] \circ
  \D{\frac{n-2}{2}} \circ \radonc \circ \sonar.
\end{equation*}
After some work, we reduce to the two-dimensional case treated in
Section \ref{s:SlantedLines}.  Effectively, our conversion of
sonar data into integrals over hyperplanes  proceeds through an
intermediate stage---integrals over cylinders.
\par
In the usual way, let $(\omega,p)$ encode a hyperplane in the
centerset of $\h{n}$. By a cylinder, with radius $r$ with axis
$(\omega,p)$, we mean any set:
\begin{displaymath}
  \left\{
    (x,y) \in \h{n} \mid
    (\omega \cdot x - p)^{2} + y^{2} = r^{2}
  \right\}.
\end{displaymath}
We encode a cylinder of radius $y$ as a triple $(\omega,p,y)$ and
write ${\cal C}[f](\omega,p,y)$ for the integral of $f$ over the
given cylinder.  (One naturally views transform ${\cal C}$ as a
hybrid of sonar and Radon.) Given  $\sonar[f]$, we can find
${\cal C}[f](\omega,p,y)$ as follows.
\begin{theorem} \label{t:SP}
For
\begin{math}
  f \in \cc{\h{n}},
\end{math}
\begin{equation} \label{e:SP}
  {\cal C}[f](\omega,p,y) =
  (\D{\frac{n-2}{2}} \circ \radonc \circ \sonar)[f](\omega,p,y)
\end{equation}
\end{theorem}
\begin{proof}
We shall actually prove the equivalent claim
$\I{\frac{n-2}{2}}\circ {\cal C}=\radonc\circ\sonar$. Combining
Equation \eqref{e:sonar} from Section \ref{s:Introduction} with
Equation \eqref{e:Radonc} from Section \ref{s:MainResult}, we
obtain an iterated integral for
\begin{math}
  ( \radonc \circ \sonar )[f](\omega,p,y)
\end{math}
in the form:
\begin{displaymath}
  \int_{\omega \cdot x = p}
  \left[
    \int_{|t| < y}
      f \left( x + t, \sqrt{y^{2} - |t|^{2}} \right) \,
      \frac{y \, dt}{\sqrt{y^{2} - |t|^{2}}}
  \right]
  \dm{x}.
\end{displaymath}
Interchanging the order of integration, which is possible because
$f$ is smooth and compactly supported, we get
\begin{displaymath}
  \int_{|t| < y}
  \left[
    \int_{\omega \cdot x = p}
      f \left( x + t, \sqrt{y^{2} - |t|^{2}} \right) \,
    \dm{x}
  \right] \,
  \frac{y \, dt}{\sqrt{y^{2} - |t|^{2}}},
\end{displaymath}
where the inside integral is a Radon transform of a shifted
function:
\begin{displaymath}
  \begin{split}
    \int_{\omega \cdot x = p}
      &f \left( x + t, \sqrt{y^{2} - |t|^{2}} \right) \,
    \dm{x} \quad
    \text{(substitute $x + t = u$)}\\
    &=
    \int_{\omega \cdot (u - t) = p}
      f \left( u, \sqrt{y^{2} - |t|^{2}} \right) \,
    \dm{u}\\
    &=
    \radonc[f]
    \left(
      \omega,p + \omega \cdot t, \sqrt{y^{2} - |t|^{2}}
    \right).
  \end{split}
\end{displaymath}
We conclude that
\begin{math}
  ( \radonc \circ \sonar )[f](\omega,p,y)
\end{math}
is the following integral over a ball
\begin{displaymath}
  \int_{|t| < y}
    \radonc[f]
    \left(
      \omega,p + \omega \cdot t, \sqrt{y^{2} - |t|^{2}}
    \right) \,
    \frac{y \, dt}{\sqrt{y^{2} - |t|^{2}}}.
\end{displaymath}
which, after switching to polar coordinates (Equation
\eqref{e:polar} from Section \ref{s:PlaneWaves}), becomes
\begin{displaymath}
  y \,
  \int_{0}^{y}
  \left[
    \int_{\theta \in \s{1}{n-2}}
      \frac{\radonc[f]
      \left(
        \omega, p + r \, (\omega \cdot \theta), \sqrt{y^{2} - r^{2}}
      \right) \, r^{n-2}}
           {\sqrt{y^{2} - r^{2}}} \,
    \ds{\theta}
  \right] \,
  dr.
\end{displaymath}
Inside the brackets, we have an integral of a plane wave over a
unit sphere.  Therefore, in light of Corollary \ref{c:John1},  the
composition
\begin{math}
  ( \radonc \circ \sonar )[f](\omega,p,y)
\end{math}
can be expressed as the following double integral:
\begin{equation} \label{e:RcS6}
  |\s{1}{n-3}| \, y \,
  \int_{0}^{y} \int_{-1}^{+1}
    \frac{
      \radonc[f]
      \left( \omega, p + r \, s, \sqrt{y^{2} - r^{2}} \right) \,
      r^{n-2}}
         {\sqrt{y^{2} - r^{2}}} \,
      \left( 1 - s^{2} \right)^{\frac{n-4}{2}} \,
  ds \, dr.
\end{equation}
The mapping
\begin{displaymath}
  (s,r) \mapsto
  \left(
    p + r \, s, \sqrt{y^{2} - r^{2}}
  \right)
\end{displaymath}
is a diffeomorphism from the rectangle
\begin{math}
  [-1,1] \times [0,y]
\end{math}
into an upper half-disk of radius $y$ centered at $p$. This
suggests the following change of variables
\begin{displaymath}
  p + r \, s = p + u, \quad
  \sqrt{y^{2} - r^{2}} = \sqrt{v^{2} - u^{2}},
\end{displaymath}
where $v \in [0,y]$ and $u < |v|$.
\par
We will now transform the integral in \eqref{e:RcS6}.  Solving for
$(s,r)$ in terms of $(u,v)$, we find that
\begin{displaymath}
  s = \frac{u}{\sqrt{y^{2} - v^{2} + u^{2}}}, \quad
  r = \sqrt{y^{2} - v^{2} + u^{2}},
\end{displaymath}
and therefore the integrand
\begin{displaymath}
  \frac{
    \radonc[f]
    \left( \omega, p + r \, s, \sqrt{y^{2} - r^{2}} \right) \, r^{n-2}}
       {\sqrt{y^{2} - r^{2}}} \,
  \left( 1 - s^{2} \right)^{\frac{n-4}{2}}
\end{displaymath}
in Equation \eqref{e:RcS6} transforms to
\begin{displaymath}
  \left( y^{2} - v^{2} \right)^{\frac{n-4}{2}} \,
  \radonc[f] \left( \omega, p + u, \sqrt{v^{2} - u^{2}} \right) \,
  \sqrt{\frac{y^{2} - v^{2} + u^{2}}{v^{2} - u^{2}}}.
\end{displaymath}
As can be found through routine differentiation, the Jacobian of
the mapping
\begin{displaymath}
  (u,v) \mapsto
  \left(
    \frac{u}{\sqrt{y^{2} - v^{2} + u^{2}}},
    \sqrt{y^{2} - v^{2} + u^{2}}
  \right).
\end{displaymath}
is the following matrix
\begin{displaymath}
  \left[%
    \begin{array}{cc}
      \frac{y^{2} - v^{2}}{\left( y^{2} - v^{2} + u^{2} \right)^{3/2}} &
      \frac{u \, v}{\left( y^{2} - v^{2} + u^{2} \right)^{3/2}} \\
      \frac{u}{\left( y^{2} - v^{2} + u^{2} \right)^{1/2}} &
      \frac{-v}{\left( y^{2} - v^{2} + u^{2} \right)^{1/2}} \\
    \end{array}%
  \right]
\end{displaymath}
whose determinant has the absolute value
\begin{displaymath}
  \frac{v}{\sqrt{y^{2} - v^{2} + u^{2}}}.
\end{displaymath}
Now \eqref{e:RcS6}, representing
\begin{math}
  ( \radonc \circ \sonar )[f](\omega,p,y),
\end{math}
takes the form
\begin{displaymath}
    |\s{1}{n-3}| \, y \,
    \int_{0}^{y}
      \left(
        y^{2} - v^{2}
      \right)^{\frac{n-4}{2}}
      \left[
        \int_{-v}^{+v}
          \radonc[f] \left( \omega, p + u, \sqrt{v^{2} - u^{2}} \right) \,
          \frac{v \, du}{\sqrt{v^{2} - u^{2}}}
      \right] \,
    dv.
\end{displaymath}
From  Definition \ref{d:I} of $\I{\nu}$ (and the surface area
formula for spheres found at the end of Section
\ref{s:PlaneWaves}) we may conclude that
\begin{displaymath}
( \radonc \circ \sonar )[f](\omega,p,y)=
    \I{\frac{n-2}{2}}
    \left\{
      \int_{-y}^{+y}
        \radonc[f] \left( \omega, p + u, \sqrt{y^{2} - u^{2}} \right) \,
        \frac{y \, du}{\sqrt{y^{2} - u^{2}}}
    \right\}.
\end{displaymath}
But now we recognize the expression in curly brackets as ${\cal
C}[f](\omega,p,y)$, as desired.
\end{proof}
To finish, note that fixing $\omega$ determines a parallel family
of cylinders.  In the definition of $\radonw[\sigma]$ we now set
$T$ equal to the set all possible $\omega$, i.e. $T=\s{1}{n-2}$.
According to Section \ref{s:SlantedLines}, the composition
\begin{math}
  \W \circ \radonw[1/y]\circ {\cal C}
\end{math}
yields the two-dimensional Radon transform of $\radonc[f]$ which
is the $n$-dimensional Radon transform $\radons[f]$.

\section{Conclusion}
\label{s:Conclusion}

The technique which proves the main theorem admits immediate
variations, if perhaps of only theoretical
interest.  For the record, we mention two.  In the sonar transform
one could replace the spheres that function as loci of integration
by other families of loci with similar scaling properties.
Alternatively, in the Radon transform, one could replace slanted planes
by cones whose axes lie in the centerset.

The authors view the methods in this paper as an expression of a
more general philosophy, under development, aimed at providing
sonar-Radon relations for more general centersets and in more
general spaces. The planar centerset case deserves an independent
treatment now because its rich structure allows for results of a
particularly explicit form and because of  potential for practical
applications.

\bibliographystyle{amsplain}
\bibliography{root}
\end{document}